\documentclass[12pt,reqno]{amsart}

\usepackage{amsmath}

\usepackage{amsfonts}

\usepackage{hyperref}

\usepackage{graphicx}

\theoremstyle{plain}

\newtheorem{defi}{Definition}

\newcommand{\brdef}{\begin{defi}}

\newcommand{\erdef}{\end{defi}}

\newcommand{\brexam}{\begin{examp}}

\newcommand{\erexam}{\end{examp}}

\newtheorem{cor}{Corollary}

\newcommand{\bcor}{\begin{cor}}

\newcommand{\ecor}{\end{cor}}

\newtheorem{thm}{Theorem}

\newcommand{\bth}{\begin{thm}}

\newcommand{\eth}{\end{thm}}

\newtheorem{lem}{Lemma}

\newcommand{\ble}{\begin{lem}}

\newcommand{\ele}{\end{lem}}

\newtheorem{pro}{Proposition}

\newcommand{\bpro}{\begin{pro}}

\newcommand{\epro}{\end{pro}}

\huge

\numberwithin{equation}{section}

\begin{document}
\title{Invariant Submanifolds of Sasakian Generalized Sasakian-space-forms}
\author{D. G. Prakasha$^{1}$, P. Veeresha$^{2, *}$, Inan Unal$^{3}$ and S. K. Hui$^{4}$}
\maketitle

\begin{center}
\pagestyle{myheadings} \markboth{D. G. Prakasha, P. Veeresha, Inan Unal and S. K. Hui}{Invariant Submanifolds of Sasakian Generalized Sasakian-space-forms}
$^{1}$Department of Mathematics, Davangere University,\\[0pt] Davangere - 577 007, India.\\[0pt]
E-mail: {\verb|prakashadg@gmail.com|, \verb|prakashadg@davangereuniversity.ac.in|}\\[0pt]
\
\\[0pt]
$^{2}$Department of Mathematics, CHRIST (Deemed to be University),\\[0pt] Bengaluru-560 029, India.\\[0pt]
E-mail: {\verb|viru0913@gmail.com|, \verb|pundikala.veeresha@christuniversity.in|}\\[0pt]
\
\\[0pt]
$^{3}$ Department of Computer Engineering,  Munzur University, \\[0pt] 
Tunceli-62000, Turkey.\\[0pt]
E-mail: {\verb|inanunal@gmail.com|}\\[0pt]
\
\\[0pt]
$^{4}$ Department of Mathematics,  The University of Burdwan, \\[0pt]
Burdwan - 713104, India.\\[0pt]
E-mail: {\verb|skhui@math.buruniv.ac.in|, \verb|shyamal_hui@yahoo.co.in|}
\end{center}
\begin{quotation}
\textbf{Abstract:}\,\  The object of this paper is to study the invariant submanifolds of Sasakian generalized Sasakian-space-form. Here, we obtain some equivalent conditions for an invariant submanifold of a Sasakian generalized Sasakian-space-forms to be totally geodesic.
\end{quotation}


\section{\bf{Introduction}}


In differential geometry, the curvature of a Riemannian manifold plays a basic role and the sectional curvature of a manifold determines the curvature tensor $\widetilde{R}$ completely. A Riemannian manifold with constant sectional curvature $c$ is called a real-space-form. In contact metric geometry, a Sasakian manifold (resp. Kenmotsu manifold) with constant $\phi$-sectional curvature is called Sasakian-space-form (resp. Kenmotsu-space-form) and it has a specific form of its curvature tensor. In 2004, Alegre et al. \cite{ABC} introduced the notion of generalized Sasakian-space-form, which can be treated as a generalization of Sasakian, Kenmotsu and Cosymplectic space-form. An almost contact metric manifold with a $\phi$-sectional curvature is called a generalized Sasakian-space-form and it is denoted by $\widetilde{M}(f_1, f_2, f_3)$. The curvature tensor of $\widetilde{M}(f_1, f_2, f_3)$ is given by \cite{ABC}:
\begin{eqnarray}\label{e:1.1}
\widetilde{R}(X, Y)Z &=& f_1 \{g(Y, Z)X-g(X, Z)Y\}\nonumber\\
&&+ f_2\{g(X, \phi Z)\phi Y-g(Y, \phi Z)\phi X+2g(X, \phi Y)\phi Z\}\\
&&+ f_3\{\eta(X)\eta(Z)Y-\eta(Y)\eta(Z)X+g(X, Z)\eta(Y)\xi -g(Y, Z)\eta(X)\xi\}\nonumber
\end{eqnarray}
for all vector fields $X, Y, Z \in \Gamma(T\widetilde{M})$, where $f_1, f_2, f_3$ are differentiable functions on $\widetilde{M}$. The $\phi$-sectional curvature of  generalized Sasakian-space-form $\widetilde{M}(f_1, f_2, f_3)$ is $f_1 +3 f_2$. Also, this notion contains a large class of almost contact manifolds. For example, any three-dimensional $(\alpha, \beta)$-trans-Sasakian manifold, where $\alpha$ and $\beta$ depending on $\xi$ is a generalized Sasakian-space-form.  In particular, if $f_1 =\frac{c+3}{4}, f_2=f_3=\frac{c-1}{4}$, then the generalized Sasakian-space-form reduces to the notion of Sasakian-space-form. The generalized Sasakian-space-forms have also been studied in \cite{AC1, AC2, AC3, HW, HP1, HP2, DGP, DCM, PN, SS, VS} and many other instances. In 2008, Alegre and carriazo \cite{AC1} studied generalized Sasakian-space-form admitting trans-Sasakian structure. In this paper, we studied generalized Sasakian-space-form admitting Sasakian structure and we call such manifold as Sasakian generalized Sasakian-space-form.
\par In modern analysis, the theory of invariant submanifolds have become today a specialized area of research due to its significant applications in applied mathematics and theoretical physics. The geometry of invariant submanifolds of almost contact manifolds was first appeared in the works of Yano and Ishihara \cite{YI}. Later, sevaral studies (see, \cite{CP, HE, HR, MTT, PVNV, TSK}) have been done on invariant submanifolds of various kinds of almost contact manifolds. For example, in \cite{MK}, Kon proved that invariant submanifolds of a Sasakian manifold are totally geodesic if their second fundamental forms are covariantly constant. De and Majhi \cite{DM} proved that an invariant submanifold $M$ of a Kenmotsu manifold is totally geodesic if and only if $Q(\sigma, R)=0$ or $Q(S, \sigma)=0$, where $\sigma$, $R$ and $S$ denote the second fundamental form, curvature tensor and Ricci tensor of $M$, respectively. Invariant submanifolds of a trans-Sasakian manifolds were studied in \cite{SS} and \cite{VS}. Further, Hu and Wang \cite{HW} investigated that an invariant submanifold of a trans-Sasakian manifold is totally geodesic if and only if $Q(S, \widetilde{\nabla}\sigma)=0$, $Q(S, \widetilde{R}\cdot \sigma)=0$, $Q(g, \widetilde{R}\cdot h)=0$, $Q(g, C\cdot \sigma)=0$ or $Q(S, C\cdot \sigma)=0$, respectively, where $C$ and $\widetilde{R}$ denote the concircular curvature tensor and curvature tensor of $\widetilde{M}$, respectively.
\par Nowadays, several authors started to study the geometry of submanifolds in the space-forms. Yildiz and Murathan \cite{YM} studied invariant submanifolds of Sasakian-space-forms. In \cite{AC2}, Alegre and Carriazo studied some submanifolds of generalized Sasakian-space-forms. Recently, Hui et. al.\cite{HUAM} studied parallel, semiparallel and 2-semiparallel invariant submanifolds of generalized Sasakian-space-forms. They also obtained the sufficient conditions of any submanifold of a generalized Sasakian-space-forms to be invariant. Motivated by the above circumstances, in this paper, we continue the study of invariant submanifolds of Sasakian generalized Sasakian-space-forms satisfying certain conditions.  
\par The paper is organized as follows: Section 2 is concerned with some preliminaries. In section 3, we study invariant submanifold of Sasakian generalized Sasakian-space-forms and prove that it is totally geodesic if and only if the second fundamental form $\sigma$ of $M$ satisfies the conditions $Q(\sigma, R)=0$, $Q(S, \sigma)=0$, $Q(S, \widetilde{\nabla}\sigma)=0$, $Q(g, \widetilde{R}\cdot \sigma)=0$, $Q(g, \mathcal{C} \cdot \sigma)=0$ or $\widetilde{R}\cdot \sigma = L_1 Q(g, \sigma)$. As a consequence of main results we obtain several corollaries. 
\section{\bf Preliminaries}
\indent \indent An odd-dimensional Riemannian manifold $(\widetilde{M}, g)$ is said to be an almost contact metric manifold \cite{DEB} if there exist on $\widetilde{M}$ a $(1, 1)$ tensor field $\phi$, a vector field $\xi$ (called the structure vector field), and a 1-form $\eta$ such that 
\begin{equation}
\eta(\xi)=1,\,\,\,\ \phi^2(X)=-X+\eta(X)\xi, \,\,\,\ \phi \xi=0,\,\,\,\ \eta \circ \phi=0,
\end{equation}
and 
\begin{equation}
g(\phi X, \phi Y)=g(X, Y)-\eta(X)\eta(Y),\,\,\,\ \eta(X)=g(X, \xi)
\end{equation}
for any vector field $X, Y \in \Gamma(TM)$.
\par On such a manifold, the fundamental 2-form $\Phi$ of $\widetilde{M}$ is defined as $\Phi(X, Y)=g(\phi X, Y)$ for any vector field $X, Y \in \Gamma(TM)$.
\par An almost contact metric manifold is called a Sasakian manifold if and only if \cite{YK}
\begin{equation}
(\widetilde{\nabla}_X \phi)Y = [g(X, Y)\xi-\eta(Y)X],\,\,\,\,\ \widetilde{\nabla}_X  \xi = -\phi X.
\end{equation}
Further, in a generalized Sasakian-space-form $\widetilde{M}(f_1, f_2, f_3)$ from (\ref{e:1.1}), we have
\begin{eqnarray}
\widetilde{S} (X, Y) &=& (2n f_1+3f_2 -f_3)g(X, Y)-\{3f_2 +(2n-1)f_3 \}\eta(X)\eta(Y),\\
\widetilde{R} (X, Y)\xi &=& (f_1 - f_3)\{\eta(Y)X-\eta(X)Y\},\\
\widetilde{R} (\xi, X)\xi &=& (f_1 - f_3)\{\eta(X)\xi - X\},\\
\widetilde{S} (\xi, \xi) &=& 2n(f_1 - f_3)
\end{eqnarray}
for all tangent vectors $X, Y \in \Gamma(TM)$, where $\widetilde{\nabla}$ denotes the covariant differentiation with respect to $g$ and $\widetilde{S}$ is the Ricci tensor of $\widetilde{M}(f_1, f_2, f_3)$.
Let $M$ be a submanifold of a $(2n+1)$-dimensional generalized Sasakian-space-form $\widetilde{M}(f_1, f_2, f_3)$. We denote by $\nabla$ and $\widetilde{\nabla}$ the Levi-Civita connection of $M$ and $\widetilde{M}$, respectively. Then, the second fundamental form $\sigma$ is given by
\begin{equation}\label{eq2.2}
\widetilde{\nabla}_X Y = \nabla_X Y + \sigma(X, Y)
\end{equation} for any vector fields $X, Y \in \Gamma(TM)$. Furthermore, for any section $N$ of normal bundle $T^{\bot}M$ we have
\begin{equation}
\label{eq2.3}
\widetilde{\nabla}_X N = -A_N X + \nabla^{\bot}_X N
\end{equation}
where $\nabla^{\bot}$ denote the normal bundle connection of $M$. If the second fundamental form $\sigma$ is identically zero then the submanifold is said to be {\it totally geodesic}. The second fundamental form $\sigma$ and $A_N$ are related by
\begin{equation}
g(\sigma(X, Y), N)=g(A_N X, Y) \nonumber
\end{equation}
for any vector fields $X,Y\in \Gamma(TM)$. For the second fundamental form $\sigma$, the first covariant derivative is given by
\begin{equation}\label{eq2.4}
(\widetilde{\nabla}_X \sigma)(Y, Z) = \nabla^{\bot}_X (\sigma(Y, Z))-\sigma(\nabla_X Y, Z)-\sigma (Y, \nabla_X Z)
\end{equation}
for any vector fields $X, Y, Z\in \Gamma(TM)$. If $\widetilde{\nabla}\sigma = 0$, then $M$ is said to have {\it parallel second fundamental form}.
\par A submanifold $M$ is said to be semi-parallel (see \cite{JD}) (resp. 2-semi-parallel, see \cite{ALMO}) if
\begin{equation}\label{eq2.6}
\widetilde{R}(X, Y)\cdot \sigma = 0\,\ (resp.\,\ \widetilde{R}(X, Y)\cdot \widetilde{\nabla}\sigma =0),
\end{equation}
holds for any vector fields $X, Y\in \Gamma(TM)$, where $\widetilde{R}(X, Y) = \widetilde{\nabla}_X \widetilde{\nabla}_Y - \widetilde{\nabla}_Y \widetilde{\nabla}_X - \widetilde{\nabla}_{[X, Y]}$ denotes the curvature tensor of the connection $\widetilde{\nabla}$. By (\ref{eq2.6}), we have
\begin{eqnarray}\label{eq2.10}
&&(\widetilde{R}(X, Y)\cdot\sigma)(U, V)\nonumber\\
&=& R^{\bot}(X, Y)\sigma(U, V) - \sigma(R(X, Y)U, V) - \sigma(U, R(X, Y)V)
\end{eqnarray}
for any vector fields $X, Y, U, V\in \Gamma(TM)$, where $R^{\bot}(X, Y) = [\nabla^{\bot}_X, \nabla^{\bot}_Y]-\nabla^{\bot}_{[X, Y]}$. Similarly, we have
\begin{eqnarray}\label{eq2.11}
&&(\widetilde{R}(X, Y)\cdot\widetilde{\nabla}\sigma)(U, V, W)\nonumber\\
 &=&R^{\bot}(X, Y)(\widetilde{\nabla}\sigma)(U, V, W)-(\widetilde{\nabla}\sigma) (R(X, Y)U, V, W) \nonumber\\
 &-& (\widetilde{\nabla}\sigma)(U, R(X, Y)V, W) - (\widetilde{\nabla}\sigma)(U, V, R(X, Y)W)
\end{eqnarray}
for any vector fields $X, Y, U, V, W\in \Gamma(TM)$, where $(\widetilde{\nabla}\sigma)(U, V, W) = (\widetilde{\nabla}_{U} \sigma)(V, W)$.
\par A submanifold $M$ is said to be pseudo-parallel (see \cite{ALM}) (resp. 2-pseudo-parallel) if 
\begin{eqnarray}
\widetilde{R}(X, Y)\cdot \sigma = L_1 Q(g, \sigma)\,\,\ (resp.\,\ \widetilde{R}(X, Y)\cdot \widetilde{\nabla}\sigma = L_1 Q(S, \widetilde{\nabla}\sigma)),\nonumber
\end{eqnarray}
holds for any vector fields $X, Y\in \Gamma(TM)$ and a smooth function $L_1$. Further, a submanbifold $M$ is said to be Ricci generalized pseudo-parallel (see \cite{MAE}) if $\widetilde{R}(X, Y)\cdot \sigma = L_1 Q(S, \sigma)$ for any vector fields $X, Y \in \Gamma(TM)$.
\par On a Riemannian manifold $M$, for a $(0, k)$-type tensor field $T$ $(k\geq 1)$ and a $(0, 2)$-type tensor field $E$, we denote by $Q(E, T)$ a $(0, k+2)$-type tensor field (see \cite{LV}), defined as follows
\begin{eqnarray}\label{eq2.12}
Q(E, T)(X_1, X_2,...,X_k;X, Y) &=& -T((X\wedge_E Y)X_1, X_2,...,X_k)\\
&-& T(X_1, (X\wedge_E Y)X_2,...,X_k)\nonumber\\
&-&...-T(X_1, X_2,...,X_{k-1}, (X\wedge_E Y)X_k),\nonumber
\end{eqnarray}
where $(X\wedge_E Y)Z = E(Y, Z)X - E(X, Z)Y$. 
\par For a $(2n+1)$-dimensional Riemannian manifold $\widetilde{M}$, the concircular curvature tensor $\mathcal{C}$ is given by
\begin{equation}\label{e:c}
\mathcal{C}(X, Y)Z = \widetilde{R}(X, Y)Z - \frac{r}{2n(2n+1)}[g(Y, Z)X - g(X, Z)Y]
\end{equation}
for any vector fields $X, Y, Z\in \Gamma(TM)$, where $r$ is the scalar curvature of $M$. We also have
\begin{equation}\label{e:c1}
(\mathcal{C}(X, Y)\cdot \sigma)(U, V) = R^{\bot}(X, Y)\sigma(U, V)-\sigma(\mathcal{C}(X, Y)U, V)-\sigma(\sigma, \mathcal{C}(X, Y)V).
\end{equation}
\section{\bf Invariant Submanifolds of Sasakian Generalized Sasakian-space-forms}
\indent\indent A submanifold $M$ of a Sasakian generalized Sasakian-space-form $\widetilde{M}(f_1, f_2, f_3)$ is called an invariant submanifold of $\widetilde{M}(f_1, f_2, f_3)$ if the characteristic vector field $\xi$ is tangent to $M$ at every point of $M$ and $\phi X$ is tangent to $M$ for any vector field $X$ tangent to $M$ at every point of $M$, that is $\phi(TM)\subset TM$ at every point of $M$.
\par Now, let $M$ be a invariant submanifold of a Sasakian generalized Sasakian-space-form $\widetilde{M}(f_1, f_2, f_3)$. Then, by Gauss formula we have
\begin{equation}\label{e:3.1}
\nabla_X \xi = -\phi X
\end{equation}
and 
\begin{equation}\label{e:3.2}
\sigma(X, \xi)=0 
\end{equation}
for any vector field $X\in\Gamma(TM)$.
Also, in an invariant submanifold $M$ of a generalized Sasakian-space-form $\widetilde{M}(f_1, f_2, f_3)$ we can deduce the following relations:
\begin{eqnarray}
\label{e:3.3} R(\xi, X)\xi &=& (f_1 -f_3)[\eta(X)\xi-X],\\
\label{e:3.4}S(X, \xi) &=& 2n(f_1 -f_3)\eta(X),\\
\label{e:3.6}\sigma(X, \phi Y) &=& \phi\,\ \sigma(X, Y)
\end{eqnarray}
for all vector fields $X, Y\in\Gamma(TM)$. Hence, we can state the following:
\begin{thm}
An invariant submanifold $M$ of a Sasakian generalized Sasakian-space-forms $\widetilde{M}$ is a Sasakian generalized Sasakian-space-form.
\end{thm}
Now, we begin with the following:
\begin{thm}\label{T2}
Let $M$ be an invariant submanifold of a Sasakian generalized Sasakian-space-form $\widetilde{M}(f_1, f_2, f_3)$ such that $f_1 \neq f_3$. Then $M$ is totally geodesic if and only if $Q(\sigma, R)=0$.
\end{thm}
{\bf Proof.} Let $M$ be an invariant submanifold of a Sasakian generalized Sasakian-space-form $\widetilde{M}(f_1, f_2, f_3)$ such that $f_1 \neq f_3$. Suppose that $Q(\sigma, R)=0$. This implies
\begin{equation}\label{e:3.7}
Q(\sigma, R)(X, Y, Z; U, V) = 0
\end{equation}
for any vector fields $X, Y, Z, U, V \in \Gamma(TM)$. By the above equation and (\ref{eq2.12}), we have
\begin{equation}\label{e:3.8}
- R((U\wedge_\sigma V)X, Y)Z- R(X, (U\wedge_\sigma V)Y)Z- R(X, Y)(U\wedge_\sigma V)Z = 0,
\end{equation}
where
\begin{equation}\label{e:3.9}
(U\wedge_\sigma V)W=\sigma(V, W)U-\sigma(U, W)V.
\end{equation}
In view of (\ref{e:3.9}), the equation (\ref{e:3.8}) can be written as
\begin{eqnarray}\label{e:3.10}
&&\sigma(U, X)R(V, Y)Z-\sigma(V, X)R(U, Y)Z+\sigma(U, Y)R(X, V)Z\nonumber\\
&&-\sigma(V, Y)R(X, U)Z+\sigma(U, Z)R(X, Y)V-\sigma(V,Z)R(X, V)U=0.
\end{eqnarray}
Setting $Z=V=\xi$ in (\ref{e:3.10}) and using (\ref{e:3.2}), we obtain
\begin{eqnarray}\label{e:3.11}
\sigma(U, X)R(\xi, Y)\xi +\sigma(U, Y)R(X, \xi)\xi =0.
\end{eqnarray}
Making use of (\ref{e:3.3}) in the above equation, we get
\begin{eqnarray}\label{e:3.12}
(f_1 - f_3)[\sigma(U, Y)X-\sigma(U, X)Y-\sigma(U, Y)\eta(X)\xi+\sigma(U, X)\eta(Y)\xi]=0.
\end{eqnarray}
Taking inner product on (\ref{e:3.12}) with $W$ and then contracting over $Y$ and $W$, we get
\begin{eqnarray}\label{e:3.13}
(2n-1)(f_1 - f_3)\sigma(U, X)=0.
\end{eqnarray}
It follows that $\sigma(U, X)=0$ for any vector fields $U, X\in \Gamma(TM)$. Hence, $M$ is totally geodesic submanifold. conversely, if $\sigma(X, Y)=0$, then from (\ref{e:3.10}), it follows that $Q(\sigma, R)=0$. This proves the theorem.
\begin{cor}
An invariant submanifold of a Sasakian-space-form is totally geodesic if and only if $Q(\sigma, R)=0$.
\end{cor}
\begin{thm}\label{T3}
Let $M$ be an invariant submanifold of a Sasakian generalized Sasakian-space-form $\widetilde{M}(f_1, f_2, f_3)$ such that $f_1 \neq f_3$. Then $M$ is totally geodesic if and only if $Q(S, \sigma)=0$.
\end{thm}
{\bf Proof.} Let us consider $Q(S, \sigma)=0$, then it follows that
\begin{eqnarray}\nonumber
Q(S, \sigma)(X, Y; U, V) = 0
\end{eqnarray}
for any vector fields $X, Y, U, V \in \Gamma(TM)$. By the above equation and (\ref{eq2.12}), we obtain
\begin{eqnarray}\label{e:3.14}
-\sigma((U\wedge_{S} V)X, Y)-\sigma(X, (U\wedge_{S} V)Y) = 0,
\end{eqnarray}
where
\begin{eqnarray}\label{e:3.15}
(U\wedge_{S} V)W=S(V, W)U-S(U, W)V.
\end{eqnarray}
Using (\ref{e:3.15}) in (\ref{e:3.14}) we have
\begin{eqnarray}\label{e:3.16}
S(U, X)\sigma(V, Y)-S(V, X)\sigma(U, Y)+S(U, Y)\sigma(X, V)-S(V, Y)\sigma(X, U)=0. \nonumber
\end{eqnarray}
Putting $V=Y=\xi$ in the above equation and then using (\ref{e:3.2}), we get
\begin{eqnarray}\label{e:3.17}
S(\xi, \xi)\sigma (U, X)=0.
\end{eqnarray}
This implies
\begin{eqnarray}\nonumber
2n(f_1 - f_3)\sigma (U, X)=0
\end{eqnarray}
for any vector fields $U, X\in\Gamma(TM)$. This completes the proof.
\begin{cor}
An invariant submanifold of a Sasakian-space-form is totally geodesic if and only if $Q(S, \sigma)=0$.
\end{cor}
\begin{thm}\label{T4}
Let $M$ be an invariant submanifold of a Sasakian generalized Sasakian-space-form $\widetilde{M}(f_1, f_2, f_3)$ such that $f_1 \neq f_3$. Then $M$ is totally geodesic if and only if $Q(S,  \widetilde{\nabla} \sigma)=0$.
\end{thm}
{\bf Proof.} Consider $Q(S, \widetilde{\nabla} \sigma)=0$, this is equivalent to 
\begin{eqnarray}\nonumber
0=Q(S, \widetilde{\nabla}_X \sigma)(W, K; U, V)
\end{eqnarray}
for any vector fields $X, W, K, U, V\in\Gamma(TM)$. We obtain directly from the above equation and (\ref{eq2.12}) that
\begin{eqnarray}\nonumber
&-&(\widetilde{\nabla}_X \sigma)(S(V, W)U, K)+(\widetilde{\nabla}_X \sigma)(S(U, W)V, K)\nonumber\\
&-&(\widetilde{\nabla}_X \sigma)(W, S(V, K)U)+(\widetilde{\nabla}_X \sigma)(W, S(U, K)V)\nonumber\\
&=& 0.
\end{eqnarray}
By the definition of $\widetilde{\nabla} \sigma$, we obtain
\begin{eqnarray}\label{e:3.18}
 && - \nabla^\perp_{X}\sigma(S(V, W)U, K)+\sigma(\nabla_X S(V, W)U, K)+\sigma(S(V, W)U, \nabla_X K)\nonumber\\
&&+ \nabla^\perp_{X}\sigma(S(U, W)V, K)-\sigma(\nabla_X S(U, W)V, K)-\sigma(S(U, W)V, \nabla_X K)\\
&& - \nabla^\perp_{X}\sigma(W, S(V, K)U)+\sigma(\nabla_X W, S(V, K)U)+\sigma(W, \nabla_X S(V, K)U)\nonumber\\
&&+ \nabla^\perp_{X}\sigma(W, S(U, K)V)-\sigma(\nabla_X W, S(U, K)V)-\sigma(W, \nabla_X S(U, K)V)\nonumber\\
&& = 0.
\end{eqnarray}
Putting $K=W=V=\xi$ in (\ref{e:3.18}) and using (\ref{e:3.2}), we obtain
\begin{eqnarray}\label{e:3.19}
S(\xi, \xi)\sigma (U, \nabla_X \xi)=0.
\end{eqnarray}
From (\ref{e:3.1}) and (\ref{e:3.4}) in (\ref{e:3.19}), we have
\begin{eqnarray}\label{e:3.20}
2n(f_1 - f_3)\sigma(U, -\phi X)=-2n(f_1 - f_3) \sigma(U, \phi X)=0.
\end{eqnarray}
Then by virtue of (\ref{e:3.6}), we have from (\ref{e:3.20}) that
\begin{eqnarray}\nonumber
2n(f_1 - f_3) \sigma(U, X)=0
\end{eqnarray}
for any vector fields $U, X\in\Gamma(TM)$. Therefore, it shows that $M$ is totally geodesic. The converse statement is trivial.
\begin{cor}
An invariant submanifold of a Sasakian-space-form is totally geodesic if and only if $Q(S, \widetilde{\nabla} \sigma)=0$.
\end{cor}
\begin{thm}\label{T5}
Let $M$ be an invariant submanifold of a Sasakian generalized Sasakian-space-form $\widetilde{M}(f_1, f_2, f_3)$ such that $f_1 \neq f_3$. Then $M$ is totally geodesic if and only if $Q(g, \widetilde{R} \cdot \sigma)=0$.
\end{thm}
{\bf Proof.} We assume that  $Q(g,  \widetilde{R} \cdot \sigma)=0$, this is equivalent to 
\begin{eqnarray}
Q(g, \widetilde{R} (X, Y)\cdot \sigma)(W, K; U, V)=0\nonumber
\end{eqnarray}
for any vector fields $X, Y, W, K, U, V\in \Gamma(TM)$. With the help of (\ref{eq2.12}), the above equation can be written as
\begin{eqnarray}
&&-(\widetilde{R}(X, Y)\cdot\sigma)(g(V, W)U, K)+(\widetilde{R}(X, Y)\cdot \sigma)(g(U, W)V, K)\nonumber\\
&&-(\widetilde{R}(X, Y)\cdot\sigma)(W, g(V, K)U)+(\widetilde{R}(X, Y)\cdot \sigma)(W, g(U, K)V)\nonumber\\
&&=0.
\end{eqnarray}
By using (\ref{eq2.11}), we obtain
\begin{eqnarray}\label{e:3.21}
&& -g(V, W)[R^\perp (X, Y)\sigma(U, K)-\sigma(R(X, Y)U, K)-\sigma(U, R(X, Y)K)]\nonumber\\
&&+ g(U, W)[R^\perp (X, Y)\sigma(V, K)-\sigma(R(X, Y)V, K)-\sigma(V, R(X, Y)K)]\nonumber\\
&&- g(V, K)[R^\perp (X, Y)\sigma(W, U)-\sigma(R(X, Y)W, U)-\sigma(W, R(X, Y)U)]\nonumber\\
&&+ g(U, K)[R^\perp (X, Y)\sigma(W, V)-\sigma(R(X, Y)W, V)-\sigma(W, R(X, Y)V)]\nonumber\\
&&=0.
\end{eqnarray}
Putting $Y=W=K=V=\xi$ in (\ref{e:3.21}) and using (\ref{e:3.2}), we obtain
\begin{eqnarray}\label{e:3.22}
\sigma(U, R(X, \xi)\xi)=0.
\end{eqnarray}
Using (\ref{e:3.3}) in (\ref{e:3.22}), we have $(f_1 - f_3)\sigma(U, X)=0$ for any vector fields $X, U\in\Gamma(TM)$. The converse is trivial.
\begin{cor}
An invariant submanifold of a Sasakian-space-forms is totally geodesic if and only if $Q(g, \widetilde{R} \cdot \sigma)=0$.
\end{cor}
\begin{thm}\label{T6}
Let $M$ be an invariant submanifold of a Sasakian generalized Sasakian-space-form $\widetilde{M}(f_1, f_2, f_3)$ such that $r\neq2n(2n+1)(f_1 - f_3) $. Then $M$ is totally geodesic if and only if $Q(g,  \mathcal{C} \cdot \sigma)=0$. 
\end{thm}
{\bf Proof.} Considering $Q(g,  \mathcal{C} \cdot \sigma)=0$, this is equivalent to 
\begin{eqnarray}\nonumber
Q(g, \mathcal{C}(X, Y)\cdot \sigma)(W, K; U, V)=0
\end{eqnarray}
for any vector fields $X, Y, W, K, U, V \in \Gamma(TM)$. Form (\ref{eq2.12}) and the above equation, we get
\begin{eqnarray}\label{e:3.23}
&& -g(V, W)(\mathcal{C}(X, Y)\cdot\sigma)(U, K)+g(U, W)(\mathcal{C}(X, Y)\cdot\sigma)(V, K)\nonumber\\
&& -g(V, K)(\mathcal{C}(X, Y)\cdot\sigma)(W, U)+g(U, K)(\mathcal{C}(X, Y)\cdot\sigma)(W, V)\nonumber\\
&& =0.
\end{eqnarray}
By the definition of $\mathcal{C}\cdot \sigma$, we obtain
\begin{eqnarray}
&& -g(V, W)[R^\perp(X, Y)\sigma(U, K)-\sigma(\mathcal{C}(X, Y)U, K)-\sigma(\mathcal{C}(X, Y)K, U)] \nonumber\\
&& +g(U, W)[R^\perp(X, Y)\sigma(V, K)-\sigma(\mathcal{C}(X, Y)V, K)-\sigma(\mathcal{C}(X, Y)K, V)] \nonumber\\
&& -g(V, K)[R^\perp(X, Y)\sigma(W, U)-\sigma(\mathcal{C}(X, Y)W, U)-\sigma(\mathcal{C}(X, Y)U, W)] \nonumber\\
&& +g(U, K)[R^\perp(X, Y)\sigma(W, V)-\sigma(\mathcal{C}(X, Y)W, V)-\sigma(\mathcal{C}(X, Y)V, W)]\nonumber\\
&&=0\nonumber.
\end{eqnarray}
Setting $Y=K=W=U=\xi$ in the above equation and then using (\ref{e:3.2}), we have
\begin{eqnarray}\label{e:3.24}
\sigma(\mathcal{C}(X, \xi)\xi, V)=0.
\end{eqnarray}
By using (\ref{e:c}), (\ref{e:3.3}) and (\ref{e:3.2}) in (\ref{e:3.24}), we immediately obtain
\begin{eqnarray}\nonumber
\left((f_1 - f_3)-\frac{r}{2n(2n+1)}\right)\sigma(X, V)=0
\end{eqnarray}
for any vector fields $X, V\in\Gamma(TM)$. This completes the proof.
\begin{cor}
An invariant submanifold of a Sasakian-space-form is totally geodesic if and only if $Q(g, \mathcal{C} \cdot \sigma)=0$, provided $r\neq 2n(2n+1)$.
\end{cor}
\begin{thm}\label{T6}
Let $M$ be an invariant submanifold of a Sasakian generalized Sasakian-space-form $\widetilde{M}(f_1, f_2, f_3)$. Then $M$ is totally geodesic if and only if the second fundamental form is pseudo-parallel, provided $L_1\neq (f_1 - f_3)$.
\end{thm}
{\bf Proof.} Let us consider a pseudo-parallel invariant submanifold of a Sasakian generalized Sasakian-space-form. Therefore, we have
\begin{eqnarray}\nonumber
\widetilde{R}\cdot \sigma=L_1 Q(g, \sigma),
\end{eqnarray}
where $L_1$ denotes the real valued function on $M$. The above equation can be written as 
\begin{eqnarray}\label{e:3.25}
(\widetilde{R}(X, Y)\cdot \sigma)(U, V)=L_1[Q(g, \sigma)(X, Y; U, V)]
\end{eqnarray}
for any vector fields $X, Y, U, V \in \Gamma(TM)$. By using (\ref{eq2.12}), we obtain from (\ref{e:3.25}) that
\begin{eqnarray}\label{e:3.26}
&&R^\perp (X, Y)\sigma(U, V)-\sigma(R(X, Y)U, V)-\sigma(U, R(X, Y)V)\\
&& =-L_1 [\sigma(X, V)g(Y, U)-\sigma(Y, V)g(X, U)+\sigma(U, X)g(Y, V)-\sigma(U, Y)g(X, V)].\nonumber
\end{eqnarray}
Putting $X=V=\xi$ in (\ref{e:3.26}), we obtain
\begin{eqnarray}\label{e:3.27}
\sigma(U, R(\xi, Y)\xi)=L_1 \sigma(U, Y).
\end{eqnarray}
Using (\ref{e:3.3}) in (\ref{e:3.27}), we get
\begin{eqnarray}\nonumber
[L_1 -(f_1 - f_3)]\sigma(U, Y)=0
\end{eqnarray}
for any vector fields $U, Y\in\Gamma(TM)$. This proves our theorem.
\begin{cor}
An invariant submanifold of a Sasakian-space-form is totally geodesic if and only if the second fundamental form is pseudo-parallel, provided $L_1 \neq 1$.
\end{cor}
For an invariant submanifold $M$ of a Sasakian generalized Sasakian-space-form $\widetilde{M}(f_1, f_2, f_3)$, from our main results we see that the following conditions are equivalent:
\begin{itemize}
\item $M$ is totally geodesic,
  \item the second fundamental form of $M$ is parallel with $f_1 \neq f_3$,
  \item the second fundamental form of $M$ is semi-parallel with $f_1 \neq f_3$,
  \item the second fundamental form of $M$ is 2-semi-parallel with $f_1 \neq f_3$,
  \item the second fundamental form of $M$ is pseudo-parallel with $f_1 - f_3 \neq L_1$,
    \item the second fundamental form of $M$ is concircularly semiparallel with $f_1 \neq f_3$ and $r\neq 2n(2n+1)(f_1 - f_3)$,
  \item the second fundamental form of $M$ is concircularly 2-semiparallel with $f_1 \neq f_3$ and $r\neq 2n(2n+1)(f_1 - f_3)$,
     \item the second fundamental form of $M$ satisfies $Q(\sigma, R)=0$ with $f_1 \neq f_3$,
    \item the second fundamental form of $M$ satisfies $Q(S, \sigma)=0$ with $f_1 \neq f_3$,
  \item the second fundamental form of $M$ satisfies $Q(S, \widetilde{\nabla}\sigma)=0$ with $f_1 \neq f_3$,
   \item the second fundamental form of $M$ satisfies $Q(g, \widetilde{R}\cdot \sigma)=0$ with $f_1 \neq f_3$,
   \item the second fundamental form of $M$ satisfies $Q(g, \mathcal{C} \cdot \sigma)=0$ with $r\neq 2n(2n+1)(f_1 - f_3)$.
\end{itemize}

%
%


\bibliographystyle{spphys}       

\begin{thebibliography}{}
%
%

\bibitem{ABC}
P. Alegre, D. E. Blair, A. Carriazo, \emph{Generalized Sasakian space forms}, Israel J. Math., 141 (1) (2004), 157-–183.

\bibitem{AC1}
P. Alegre, A. Carriazo, \emph{Structures on generalized Sasakian space forms}, Differential Geom. Appl., 26 (6) (2008), 656-–666.

\bibitem{AC2}
P. Alegre, A. Carriazo,  \emph{Submanifolds of generalized Sasakian space forms}, Taiwanese J. Math., 13 (3) (2009), 923–-941.

\bibitem{AC3}
P. Alegre, A. Carriazo,  \emph{Generalized Sasakian space forms and conformal change of metric}, Results Math., 59 (3-4) (2011), 485–-493.

\bibitem{ALMO}
K. Arslan, U. Lumiste, C. Murathan, C. Ozgur, \emph{2-semiparallel surfaces in space forms, I: two particular cases}, Proc. Estonian Acad. Sci. Phy. Math., 49 (3) (2000), 139--148.

\bibitem{ALM}
A. C. Asperti, G. A. Lobos, F. Mercuri, \emph{Pseudo-parallel immersions in space forms}, Mat. Contemp., 17 (1999) 59--70.

\bibitem{DEB}
D. E. Blair, Contact manifolds in Riemannian geometry, Lecture Notes in Mathematics, Vol. 509, Springer-Verlag, 1976.

\bibitem{CP}
D. Chinea, P. S. Prestelo, \emph{Invariant submanifolds of a trans-Sasakian manifold}, Publ. Math. Debrecen, 38 (1-2) (1991), 103--109.

\bibitem{DM}
U. C. De, P. Majhi, \emph{On invariant submanifolds of Kenmotsu manifolda}, J. Geom., 106 (1) (2015), 109--122.


\bibitem{JD}
J. Deprez, \emph{Semi-parallel surfaces in Euclidean space}, J. Geom., 25 (2) (1985), 192--200.

\bibitem{HE}
H. Endo, \emph{Invariant submanifolds in a contact metric manifolds}, Tensor, 43 (1986), 83--87.

\bibitem{HW}
C. Hu, Y. Wang, \emph{A note on invariant submanifolds of trans-Sasakian manifolds}, Int. Electron. J. Geom., 9 (2) (2016), 27--35.

\bibitem{HP1}
S. K. Hui, D. G. Prakasha, \emph{On the C-Bochner curvature tensor of generalized Sasakian-space-forms}, Proc. Natl. Acad. Sci., India, Sect. A Phys. Sci., 85 (3) (2015), 401--405.

\bibitem{HP2}
S. K. Hui, D. G. Prakasha, V. Chavan \emph{On generalized $\phi$-recurrent generalized Sasakian-space-forms}, Thai J. Math., 15 (2) (2017), 323--332.

\bibitem{HR}
S. K. Hui, J. Roy, \emph{Invariant and anti-invariant submanifolds of special quasi-Sasakian manifolds}, J. Geom., 109 (37) (2018), 1-16. https://doi.org/10.1007/s00022-018-0442-2.


\bibitem{HUAM}
S. K. Hui, S. Uddin, A. H. Alkhaldi, P. Mandal, \emph{Invariant submanifolds of generalized Sasakian-space-forms}, Int. J. Geom. Methods Mod. Phys., 15 (9) (2018). DOI: 10.1142/S0219887818501499.


\bibitem{MK}
M. Kon, \emph{Invariant submanifolds of normal contact metric manifolds}, Kodai Math. Sem. Rep., 25 (1973), 330--336.

\bibitem{MTT}
B. C. Montano, L. D. Terlizzi, M. M. Tripathi, \emph{Invariant submanifolds of Contact $(\kappa, \mu)$-manifolds}, Glasgow Math. J., 50 (2008), 499--507.

\bibitem{MAE}
C. Murathan, K. Arslan, E. Ezentas, \emph{Ricci generalized pseudo-symmetric manifolds}, Differ. Geom. Appl., 99--108, Matfyzpress, Prague, 2005.

\bibitem{DGP}
D. G. Prakasha, \emph{On generalized Sasakian space forms with Weyl-conformal curvature tensor}, Lobachevskii J. Math., 33 (3) (2012), 223-–228.

\bibitem{DCM}
D. G. Prakasha, V. Chavan, K. Mirji, \emph{On the $W_5$-curvature tensor of generalized Sasakian-space-forms}, Konuralp J. Math., 4 (1) (2016), 45-–53.

\bibitem{PN}
D. G. Prakasha, H. G. Nagaraja, \emph{On quasi-conformally flat and quasi-conformally semisymmetric generalized Sasakian space forms}, CUBO, 15 (3) (2013), 59-–70.

\bibitem{PVNV}
D. G. Prakasha, A. T. Vanli, M. Nagaraja, P. Veeresha, \emph{Invariant submanifolds of $(\epsilon)-$Sasakian manifolds}, Ann. Univ. Craiova Math. Comput. Sci. Ser., 47 (2) (2020), 346-–356.


\bibitem{SS}
A. Sarkar, M. Sen, \emph{On invariant submanifolds of trans-Sasakian manifolds}, Proc. Est. Aacad. Sci., 61 (1) (2012), 29--37.

\bibitem{TSK}
M. M. Tripathi, T. Sasahara, J. S. Kim, \emph{On invariant submanifolds of contact metric manifolds}, Tsukuba J. Math., 29 (2) (2005), 495--510.

\bibitem{VS}
A. T. Vanli, R. Sari, \emph{Invariant submanifolds of trans-Sasakian manifolds}, Differ. Geom. Dyn. Syst., 12 (2010), 227--288.

\bibitem{LV}
L. Verstraelen, \emph{Comments on pseudosymmetry in the sense of Ryszard Deszcz}, Geometry and Topology of submanifolds, 6 (1) (1994), 199--209.

\bibitem{YI}
K. Yano, S. Ishihara, \emph{Invariant submanifolds of an almost contact manifold}, Kodai Math. Sem. Rep., 21 (1969), 350--364.

\bibitem{YK}
K. Yano, M. Kon, \emph{Structures on manifolds}, World Scientific, (1985).

\bibitem{YM}
A. Yildiz, C. Murathan, \emph{Invariant submanifolds of Sasakian-space-forms}, J. Geom., 95 (1-2) (2009), 135--150.


\end{thebibliography}


\end{document}